\documentclass{amsart}

\usepackage{fullpage}

\usepackage{amsmath,amsfonts,amssymb}

\usepackage{graphicx}                          
\usepackage{setspace}  
\usepackage{caption}
\usepackage{listings}
\usepackage{tabularx}
\usepackage{xcolor}
\usepackage{indentfirst}
\usepackage{cite}
\usepackage{subcaption}

\usepackage{overpic}

\usepackage{bm}

\newtheorem{remark}{Remark}[section]

\newtheorem{theorem}{Theorem}[section]

\usepackage{mhchem}

\newcommand\dd{\mathrm{d}}
\newcommand\pp{\partial}
\newcommand\x{\bm{x}}
\newcommand\uvec{\mathbf{u}}

\begin{document}

\title{On a Continuum Model for Random Genetic Drift:  A Dynamic Boundary Condition Approach}

\author{Chun Liu}
\address{Department of Applied Mathematics, Illinois Institute of Technology, Chicago, IL.}
\email{cliu124@iit.edu}

\author{Jan-Eric Sulzbach}
\address{Department of Mathematics, Technical University Munich, Munich, Germany.}
\email{janeric.sulzbach@ma.tum.de}

\author{Yiwei Wang}
\address{Department of Mathematics, University of California, Riverside, CA.}
\email{yiweiw@ucr.edu}



\begin{abstract}
  We propose a new continuum model for random genetic drift by employing a dynamic boundary condition approach. The model can be viewed as a regularized version of the Kimura equation and admits a continuous solution. We establish the existence and uniqueness of a strong solution to the regularized system. Numerical experiments illustrate that, for sufficiently small regularization parameters, the model can capture key phenomena of the original Kimura equation, such as gene fixation and conservation of the first moment.
\end{abstract}

\maketitle

\section{Introduction}

Genetic drift is a fundamental process in molecular evolution \cite{ewens2004mathematical}. It refers to the random fluctuations in allele frequencies within a population over time. Typical mathematical models for genetic drift include the Wright--Fisher model \cite{fisher1923xxi, wright1931evolution} and its diffusion limit, known as the Kimura equation \cite{kimura1962probability, wright1945differential}.

Without mutation, migration, and selection process, the Kimura equation is a Fokker-Planck type equation, given by \cite{ewens2004mathematical}
\begin{equation}\label{Kimura_eq}
  \pp_t \rho = \pp_{xx} ( x(1 - x) \rho) \ , \quad x \in (0, 1) \ , \  t >0 \ .
 \end{equation}
Here, the variable $x$ is the fraction of the focal allele $A_1$ in the population and $(1 - x)$ is that of $A_2$, the function $\rho(x, t)$ is the probability of finding a relative composition $x \in (0,1)$ of gene $A_1$ at time $t$. 
Although (\ref{Kimura_eq}) is a linear PDE of $\rho(x, t)$, due to the degeneracy of $x(1 - x)$ at the boundary $x = 0$ and $x = 1$ \cite{epstein2010wright, epstein2013degenerate}, even the appropriate boundary conditions and solution space of the Kimura equation are unclear \cite{feller1951diffusion}. In the 1950s, Kolmogorov suggested that equation \eqref{Kimura_eq} is only reasonable for $x$ not too close to $0$ and $1$ \cite{Shiryayev1992}.

To maintain the biological significance, $\rho(x, t)$ should satisfy the mass conservation
\begin{equation}\label{conservation_law_1}
   \int_{0}^1 \rho(x, t) \dd x = 1 \ ,
 \end{equation}
 which suggests a non-flux boundary to (\ref{Kimura_eq}), given by
 \begin{equation}\label{BC_non_flux}
  \pp_x(x(1 - x) \rho) = 0 \ , \quad x = 0 ~\text{or}~ 1 \ , \quad \forall t \ .
 \end{equation}
 Moreover, the biologically relevant solution of the Kimura equation (\ref{Kimura_eq}) needs to satisfy an additional conservation law: 
\begin{equation}\label{conservation_law_2}
\frac{\dd}{\dd t} \int_{0}^1 x \rho(x, t) \dd x = 0\ ,
\end{equation}
known as the conservation of fixation probability in biology. Here, $\psi(x) = x$ is the fixation probability function, which describes the probability of allele $A_1$ fixing in a population while allele $A_2$ goes extinct, under the condition of starting from an initial composition of $x$. It can be noticed that in the pure drift case, (\ref{conservation_law_2}) is the conservation of the first moment of $\rho(x, t)$.
Formally, the conservation of the fixation probability can be obtained through integration by parts with the non-flux condition (\ref{BC_non_flux})
\begin{equation}
  \begin{aligned}
  & \frac{\dd}{\dd t} \int_0^1 x \rho(x, t) \dd x = \int_0^1 x \pp_{xx} ( x(1 - x) \rho)  \dd x = \int_{0}^1 - \pp_x (x (1 - x) \rho) \dd x  + x \pp_x (x (1 - x) \rho) \Big|_{0}^1 = 0 \\ 
  \end{aligned}
\end{equation}
 
However, the boundary condition (\ref{BC_non_flux}) leads to a finite-time blow-up of the solution of (\ref{Kimura_eq}) on the boundary \cite{chalub2009non}. In \cite{chalub2009non, mckane2007singular}, the authors show that for a given $\rho_0 \in \mathcal{BM}^+([0, 1])$,
 there exists a unique measure-valued weak solution to (\ref{Kimura_eq}) with $\rho(x, t) \in L^{\infty} ([0, \infty), \mathcal{BM}^+([0, 1]))$ that satisfies two conservation laws (\ref{conservation_law_1}) and (\ref{conservation_law_2}), and the solution $\rho(x, t)$ can be expressed as
\begin{equation}\label{measure_solution}
 \rho(x, t) = q(x, t) + a(t) \delta_0  + b(t) \delta_1 \ .
 \end{equation}
 Here, $\mathcal{BM}^+([0, 1])$ is the space of all (positive) Radon measures on $[0,1]$, $\delta_0$ and  $\delta_1$ are Dirac delta functions at 0 and 1 respectively, and $q(x, t) \in C^{\infty} (\mathbb{R}^+; C^{\infty}([0, 1]))$ is a classical solution to (\ref{Kimura_eq}) without boundary conditions. Moreover, it is proved that \cite{chalub2009non}, as $t \rightarrow \infty$, $q (x, t) \rightarrow 0$ uniformly, and 
 $a(t)$ and $b(t)$ are monotonically increasing functions such that 
  \begin{equation}\label{eq_a_b}
 \begin{aligned}
   & a^{\infty} = \lim_{t \rightarrow \infty } a(t) = \int_{0}^1 (1 - x) \rho_0(x) \dd x\ , \quad b^{\infty} = \lim_{t \rightarrow \infty } b(t) = \int_{0}^1 x \rho_0(x) \dd x \  . \\
 \end{aligned}
 \end{equation}

The measure-valued solution (\ref{measure_solution}) impose a difficulty to study the Kimura equation both numerically \cite{duan2019numerical} and theoretically \cite{casteras2022hidden}. For instance, although a lot of numerical methods for the Kimura equation have been developed, ranging from finite volume method \cite{xu2019behavior, zhao2013complete}, finite Lagrangian methods \cite{duan2019numerical}, an optimal mass transportation method \cite{carrillo2022optimal}, and SDE--based simulation methods \cite{dangerfield2012boundary, jenkins2017exact},
it is difficult to get a good numerical approximation to (\ref{measure_solution}), which includes the Dirac delta function, and accurately capture the dynamics of $a(t)$ and $b(t)$. 

The purpose of this paper is to propose a new continuum model for random genetic drift by incorporating a dynamic boundary approach. The key idea is to only consider the Kimura equation on $(\delta, 1 - \delta)$ for a given small $\delta > 0$, and use
a dynamic boundary condition to describe the jump process between bulk and surface, i.e., the flux between $\rho(x, t)$ and the the boundary states $a(t)$ and $b(t)$.
The new model is given by
\begin{equation} \label{eq_reg_kimura_intro}
  \begin{cases}
    & \rho_t = \pp_{xx} (x (1 - x) \rho), \quad x \in (\delta, 1 - \delta), \\ %
    & \pp_x (x (1 - x) \rho) \Big|_{x = \delta} = a'(t), \quad \pp_x (x (1 - x) \rho) \Big|_{x = 1 - \delta} = - b'(t) \\ 
    & a'(t) = -  ( (\epsilon a) - \rho(\delta , t) )  \\
    & b'(t) = -  ( (\epsilon b) - \rho(1 - \delta , t) ).  \\
  \end{cases}
 \end{equation}
Here, $a(t)$ and $b(t)$ are probability at $x = 0$ and $x = 1$ respectively, $\epsilon > 0$ is an additional small parameter that is introduced such that the overall system has a variational structure.
Such a regularized Kimura equation admits a strong global solution for fixed $\epsilon$ and $\delta$.  Although it is difficult to prove theoretically, numerical tests show that the qualitative behaviors of the original Kimura equation can be well captured with small $\epsilon$ and $\delta$. Specifically, numerical results show that $\rho(x, t)$ serves as an approximation to $q(x, t)$ in (\ref{measure_solution}), while $a(t)$ and $b(t)$ capture the singular behavior of the Dirac-delta function at the boundary in (\ref{measure_solution}).

 The rest of this paper is organized as follows. In Section 2, we introduce some background, including the Wright-Fisher model, the variational structure of Kimura equation, and the dynamic boundary approach for generalized diffusions. The new continuum model is presented in Section 3. The existence and uniqueness of the strong solution of the regularized system is shown in Section 4. Finally, we perform a numerical study on the regularized system, and demonstrate the effects of $\delta$ and $\epsilon$, as well as the ability of the new model in capturing the key feature of the original Kimura equation.

\section{Background}

\subsection{From a Wright-Fisher model to the Kimura equation}
We first briefly review the formal derivation of the Kimura equation from the Wright-Fisher model.
Consider two competing alleles, denoted by $A_1$ and $A_2$, in a diploid population with fixed size $N$ (i.e., total $2N$ alleles).
Assume $u$ and $v$ are the probabilities of mutations $A_1 \rightarrow A_2$ and $A_2 \rightarrow A_1$ respectively.
Let $X_k$ be the portion of individuals of type $A_1$ in generation $k$.  The Wright-Fisher model describes 
random fluctuations in the genetic expression by a discrete-time, discrete-state Markov chain. The discrete state space is defined as $S = \{0, \frac{1}{2N}, \ldots, 1\}$, and the transition probability between two states is given by \cite{feller1951diffusion, ethier1977error}
\begin{equation}
\mathbb{P} \left( X_{k+1} = \frac{j}{2N}   ~\Big|  X_k = \frac{i}{2 N} \right) = \tbinom{2N}{j} p_i^j (1 - p_i)^{2N - j}\ , 
\end{equation}
where
\begin{equation}
    p_i = (1 - u) \frac{i}{2 N} + v \left(1 - \frac{i}{2 N} \right)\ .
\end{equation}
In the case that $u = v  = 0$ (without mutations), we have
\begin{equation}\label{Kimura_0}
  \mathbb{P} \left( X_{k+1} = \frac{j}{M}   ~\Big|  X_k = \frac{i}{2 N} \right) = \tbinom{2N}{j} \left(\frac{i}{2N} \right)^j \left(1 - \frac{i}{2N} \right)^{2N - j}\ ,
\end{equation}
which is known as the pure drift case. Notice that in the pure drift case, the jump rate becomes 0 if $X_k = 0$ or $X_k = 1$, which are two absorbing states in the system.  These states represent allele fixation, where one allele becomes fixed (frequency 1) and the other lost (frequency 0). The system does not permit fixation if $u > 0$ and $v > 0$.

When the population size $N$ is large, the Wright-Fisher model can be approximated by a continuous--state, continuous--time process $X_t$, which represents the proportion of alleles of type $A_1$.
The dynamics of $X_t$ can be described by an SDE \cite{dangerfield2012boundary}
\begin{equation}
\dd X_t = (u - (u + v) X_t) \dd t + \sqrt{X_t (1 - X_t)} \dd W_t\ ,
\end{equation}
where $\dd W_t$ is the standard Brownian motion. The corresponding Fokker-Planck equation is given by \cite{ethier1977error, tran2013introduction}
\begin{equation}
 \rho_t +  \pp_x ((u(1 - x) - v x ) \rho) = \frac{1}{2} \pp_{xx} ( x (1 - x) \rho)\ , \quad x \in (0, 1), \quad t >0,
\end{equation}
where $\rho(x, t)$ is the probability density of $X_t$.
For the pure drift case, we obtain a Kimura equation
\begin{equation}\label{Kimura_1}
  \pp_t \rho = \frac{1}{2} \pp_{xx} ( x(1 - x) \rho), \quad x \in (0, 1), \quad t >0,
 \end{equation}

It is convenient to re-scale the time by letting $t' = 2 t$, and the Kimura equation \eqref{Kimura_1} becomes \eqref{Kimura_eq}.

\subsection{Formal Variational Structure}
Formally, the Kimura equation (\ref{Kimura_eq}) can be viewed as a generalized diffusion equation, derived from the energy-dissipation law \cite{duan2019numerical}
\begin{equation}\label{Energy-Dissipation}
 \frac{\dd}{\dd t} \int_{0}^1 \rho \ln \left( x(1 - x) \rho \right) \dd x = - \int_{0}^1 \frac{\rho}{x(1 - x)}  |u|^2   \dd x\ ,
\end{equation}
where the distribution $\rho$ satisfies the kinematics
\begin{equation}\label{kinemtics}
 \begin{aligned}
   & \pp_t \rho + \pp_x (\rho u) = 0\ , \quad x \in (0, 1) \\
   & \rho u  = 0, \quad x = 0~\text{or}~1 \ . \\
 \end{aligned}
\end{equation}
and $u$ is the macroscopic deterministic velocity for the diffusion.
The non-flux boundary condition $\rho u = 0$ at $x = 0$ and $1$ guarantees the mass conservation.

To obtain the Kimura equation from the energy-dissipation law (\ref{Energy-Dissipation}), one needs to introduce a flow map denoted as $x(X, t): [0,1] \rightarrow [0, 1]$ associated with the velocity field $u(x, t)$. For given $X$, the flow map satisfies the following ordinary differential equation:
\begin{equation}
\frac{\dd}{\dd t} x(X, t) = u(x(X, t), t), \quad x(X, 0) = X,
\end{equation}
where $X$ represents the Lagrangian coordinate, and $x$ represents the Eulerian coordinate. Due to the mass conservation, $\rho(x(X, t), t)$ is determined by the flow map $x(X, t)$ through
\begin{equation}
\rho(x(X, t), t) = \rho_0(X) / \det F(X, t)\ , \quad \det F(X, t) = \pp_X x(X, t)\ , \quad X \in (0, 1)\ ,
\end{equation}
where $\rho_0(X)$ is the initial density.
Recall $u(x(X, t), t) = x_t(X, t)$, the energy-dissipation law can be written as
\begin{equation}\label{ED_Lagrangian}
\frac{\dd}{\dd t} \int \rho_0 \ln(x(1 - x)) + \rho_0 \ln (\rho_0 / \pp_X x)  \dd X = - \int_{0}^1 \frac{\rho_0}{x(1 - x)} |x_t|^2 \dd x
\end{equation}
in Lagrangian coordinates, which can be interpreted as an $L^2$-gradient flow in terms of $x(X, t)$.
By a standard energetic variation procedure (see \cite{duan2019numerical}), we can derive the force balance equation
\begin{equation}
\frac{\rho}{x(1 - x)} u = - \rho  \pp_x( \ln  (\rho x(1 - x)))\ ,
\end{equation}
which can be simplified as
\begin{equation}\label{velocity_a}
  \rho u  = -  \pp_x (x (1 - x)p)\ . 
\end{equation}
Combining the velocity equation with the kinematics (\ref{kinemtics}), we can recover the original Kimura equation.
The variational structure (\ref{ED_Lagrangian}) naturally leads to the Lagrangian algorithms for the Kimura equation developed in \cite{duan2019numerical,carrillo2022optimal}.

Alternatively, (\ref{Energy-Dissipation}) can be interpreted a Wasserstein type of gradient flow, with the transport distance defined by 
\begin{equation}\label{WS_distance}
  \begin{aligned}
 & d^2(\rho_1, \rho_2) = \mathop{\min}_{(\rho, \uvec)} \int_{0}^1 \int_0^1 \frac{\rho}{x(1 - x)} |u|^2 \dd x \dd t~,  \\
 & \text{subject to}~ \rho_t + \pp_x (\rho u) = 0~,\quad \rho(x, 0) = \rho_1~, \quad \rho(x, 1) = \rho_2. \\
  \end{aligned}
\end{equation}
The distance (\ref{WS_distance}) is known as Wasserstein-Shahshahani distance \cite{carrillo2022optimal, chalub2021gradient, casteras2022hidden}.

As pointed out in \cite{casteras2022hidden}, the variational structure presented in this subsection is rather formal. Although we can define $U(x) = \ln (x (1 - x))$ \cite{casteras2022hidden, tran2015free}, which plays the role of internal energy as in the standard Fokker-Planck equation, the system doesn't admit a unique equilibrium distribution $\rho^{\rm eq} \propto \exp( - U)$, as $\int_{0}^1 \frac{1}{x (1 - x)} \dd x$ is unbounded.

\subsection{Dynamic boundary condition}

As suggested in \cite{casteras2022hidden}, a natural attempt at compensating for singularities is to relax the boundary condition by taking into account of the bulk/surface interaction. This leads to the use of dynamic boundary conditions \cite{knopf2021phase, monsaingeon2021new}, in which the boundary values evolve according to their own differential equations, often coupled with the interior dynamics.

Before we present the application of the dynamic boundary condition approach to the Kimura model, we first briefly review the dynamic boundary approach for generalized diffusions in this subsection. Consider a bounded domain $\Omega$, and let $\pp \Omega$ be the boundary of $\Omega$.
Classical PDE models on $\Omega$ often impose Dirichlet, Neumann, and Robin boundary conditions for the physical variable $\rho$. The basic fundamental assumption behind this is that $\rho$ is regular enough (e.g. $\rho \in C^m (\bar{\Omega})$ for some $m$), allowing $\rho|_{\pp \Omega}$ to be defined as the trace of $\rho$, and the boundary condition to be specified in terms of $\rho|_{\Gamma}$. However, as in the classical Kimura equation, it may not be always possible to take the trace.

The idea of dynamic boundary condition is to introduce another function $\sigma \in C(\pp \Omega)$ to describe the surface densities on the boundary, and view the exchange between bulk and surface density as a chemical reaction $\rho \ce{<=>} \sigma$ \cite{knopf2021phase, wang2022some}. 
Due to the mass conservation, $\rho$ and $\sigma$ satisfies
\begin{equation}
\frac{\dd}{\dd t} \left(\int_{\Omega} \rho \dd \x + \int_{\pp \Omega} \sigma \dd S \right)  = 0,
\end{equation}
which leads to the kinematics in Eulerian coordinates
\begin{equation}
  \begin{aligned}
    & \rho_t + \nabla \cdot (\rho \uvec) = 0, ~~ \x \in \Omega, \\
    & \rho \uvec \cdot {\bm \nu} = \dot{R}, \quad \sigma_t + \nabla_{\pp \Omega} \cdot (\sigma {\bm v}) = \dot{R}, ~~ \x \in \pp \Omega \\
  \end{aligned}
\end{equation}
where ${\bm \nu}$ is the outer normal of $\Omega$ and $R$ is the reaction trajectory for the chemical reaction $\rho \ce{<=>} \sigma$. The reaction trajectory $R$, representing the number of reactions in the forward direction that have occurred by time $t$, is analogous to the flow map in mechanical systems \cite{oster1974chemical, wang2020field}. The quantity $\dot{R}$ corresponds to the reaction velocity or reaction rate \cite{wang2020field}.

In general, systems with dynamic boundary condition can be modeled through an energy-dissipation law
\begin{equation}
  \frac{\dd}{\dd t} \left( \mathcal{F}_b(\rho) +  \mathcal{F}_s(\sigma) \right) = -  \left( \int_{\Omega} \eta_b (\rho) |\uvec|^2 \dd \x + \int_{\pp \Omega} \eta_s (\sigma) |{\bm v}|^2  + \dot{R} \Gamma(R, \dot{R}) \dd S \right)\ .
\end{equation}
Here, $\mathcal{F}_b(\rho) $ and $\mathcal{F}_s(\sigma)$ are free energies in the bulk and surface respectively, $\eta_b (\rho) > 0$ and $\eta_{s} (\sigma) > 0$ are friction coefficients for bulk and surface diffusions, $\dot{R} \Gamma(R, \dot{R}) \geq 0$ is the dissipation due to the bulk/surface interaction, which in general is non-quadratic in terms of $\dot{R}$ \cite{wang2020field}. The variational procedures lead to the force balance equations for the mechanical and chemical parts (see \cite{wang2020field} for details)
\begin{equation}
  \begin{cases}
    & \eta_{b} (\rho) \uvec = - \rho \nabla \mu_b (\rho), \quad \mu_b = \frac{\delta \mathcal{F}_b}{\delta \rho} \quad \x \in \Omega \\
    & \eta_s (\sigma) {\bm v} = - \sigma \nabla_{\Gamma} \mu_s (\sigma), \quad \mu_s = \frac{\delta \mathcal{F}_s}{\delta \sigma}, \quad \x \in \Gamma \\
    & \Gamma(R, \dot{R}) = - (\mu_s(\sigma) - \mu_b(\rho)), \quad \x \in \Gamma  \\
  \end{cases}
\end{equation}
where $\mu_s(\sigma) - \mu_b(\rho)$ is the affinity of the bulk-surface reaction. Different choices of free energy and dissipation lead to different systems \cite{wang2022some}. 


\section{A regularized Kimura equation}
In this section, we propose a regularized Kimura equation by applying a dynamic boundary condition approach in this specific one dimensional setting.

\subsection{Model derivation}

The key idea is to only consider the Kimura equation on $(\delta, 1 - \delta)$ for a given small $\delta > 0$, and use a dynamic boundary condition to describe the fixation on the boundary. Specially, for a given $\delta > 0$, let $\rho(x, t)$ represent the probability that the gene frequency is equal to  $x \in (\delta, 1 - \delta)$. We denote the probabilities of the gene frequency being $x = 0$ and $x = 1$ as $a(t)$ and $b(t)$, respectively.
Due to the conservation of mass, we have
\begin{equation}\label{conservation_of_mass}
 \begin{aligned}
   \frac{\dd}{\dd t}  \left(  \int_{\delta}^{1-\delta} \rho(x, t) \dd x + a(t) + b(t) \right) = 0
 \end{aligned}
\end{equation}
which leads to the following kinematics
\begin{equation}\label{kinemtics_DB}
 \begin{aligned}
  & \pp_t \rho + \pp_x (\rho u) = 0, \quad x \in (\delta, 1 - \delta) \\
  &  \rho u(\delta, t) = - \dot{R}_0(t), \quad \rho u(1 - \delta, t) = \dot{R}_1 (t) \\
  &  a'(t) = \dot{R}_0(t), \quad    b'(t)  = \dot{R}_1(t) \\
 \end{aligned}
\end{equation}
Here, $R_0(t)$ and $R_1(t)$ denotes the reaction trajectory for the bulk-boundary interactions $x = \delta$ with $x = 0$ and $x = 1 - \delta$ with $x = 1$, i.e., $\rho(\delta) \ce{<=>} a(t)$ and $\rho(1 - \delta) \ce{<=>} b(t)$. 

The overall system can be modeled through the energy-dissipation law
\begin{equation}\label{ED_RK}
\begin{aligned}
 & \frac{\dd}{\dd t} \int_{\delta}^{1 - \delta} \rho \ln \left( x(1 - x) \rho \right)\dd x + G_0(a) + G_1(b)  = - \int_{\delta}^{1 - \delta} \frac{\rho}{x(1 - x)}  |u|^2   \dd x - \Gamma_0(R_0, \dot{R_0}) - \Gamma_1(R_1, \dot{R_1})
 \end{aligned}
\end{equation}
where $G_0(a)$ and $G_1(b)$ are free energy on the boundary, and $\Gamma_0(R_0, \dot{R}_0)$ and $\Gamma_1(R_1, \dot{R}_1)$ are dissipations on due to the jump between bulk and surface. The energy-dissipation law in the bulk region ($x \in (\delta, 1- \delta)$) is exactly the same as that in the original Kimura equation (\ref{Energy-Dissipation}).
The remaining question is how to choose $G_i (i = 0, 1)$ and $\Gamma_i (i = 0, 1)$ in order to capture the qualitative behavior of the original Kimura equation. In the current work, we take
\begin{equation}\label{Form_G}
 G_0(q) = G_1(q) = G(q) =  q \ln ( \eta(\epsilon) \delta (1 - \delta) q   )\ , 
\end{equation}
and
\begin{equation}
 \Gamma_0(R, \dot{R})  = \dot{R} \ln \left( \frac{\dot{R}_t}{\kappa(\epsilon) a}  + 1  \right), \quad  \Gamma_1(R, \dot{R})  = \dot{R} \ln \left( \frac{\dot{R}}{\kappa(\epsilon) b}  + 1  \right) \ .
\end{equation}
Here, $\eta(\epsilon) > 0$ determines the boundary-surface equilibrium , $\kappa(\epsilon)$ is the reaction rate from surface to bulk. Moreover, we take $\kappa(\epsilon) = \epsilon$ and $\eta(\epsilon) = \frac{\epsilon}{\gamma}$ with $\gamma > 0$ being a constant. We'll explore other choices in future work.

By an energetic variational procedure \cite{wang2022some}, we can obtain the velocity equation
\begin{equation}\label{velocity}
 \rho u = -  \pp_x (x (1 - x) \rho), \quad x \in (\delta, 1 - \delta), 
\end{equation}
and the equations for reaction rates
\begin{equation}\label{Eq_R}
 \begin{aligned}
   & \ln \left( \frac{\dot{R}_0}{\epsilon a} + 1   \right) =  - (\ln ( \tfrac{\epsilon}{\gamma} a) -  \ln \rho(\delta, t))  \\
   & \ln \left( \frac{\dot{R}_1}{\epsilon b}  + 1  \right) = - (\ln ( \tfrac{\epsilon}{\gamma} b) - \ln \rho(1 - \delta, t) )\ .  \\
 \end{aligned}
\end{equation}
One can rewrite (\ref{Eq_R}) as
\begin{equation}\label{Eq_RR}
\dot{R}_0 = \gamma \rho(\delta, t) - \epsilon a, \quad  \dot{R}_1 =  \gamma \rho(1 - \delta, t) - \epsilon b \ .
\end{equation}
Combing (\ref{velocity}) and (\ref{Eq_RR}) with the kinematics (\ref{kinemtics_DB}), one arrives the final equation 
\begin{equation}\label{Eq_Final_LMA}
 \begin{cases}
   &  \pp_t \rho = - \pp_x (\rho u), \quad   x \in [\delta, 1 - \delta] \\
   & \rho u = -  \pp_x (x (1 - x) \rho), \quad x \in (\delta, 1 - \delta), \\ %
   & \rho u(\delta, t) = - a'(t), \quad \rho u(1 - \delta, t) = b'(t) \\
   & a'(t) =  \gamma \rho(\delta, t) - \epsilon a  \\
   & b'(t) =  \gamma  \rho(1 - \delta, t) - \epsilon b  \\
 \end{cases}
\end{equation}

For the regularized system (\ref{Eq_Final_LMA}), since
\begin{equation}
  a'(t) + \epsilon a =  \gamma \rho(\delta, t)
\end{equation}
we have
\begin{equation}
  a(t) =   e^{-\epsilon t} \int_{0}^t \gamma e^{\epsilon s} \rho(\delta, s) \dd s  
\end{equation}
if the initial condition $a(0) = 0$.
As a consequence,
\begin{equation}
a'(t) = -\epsilon e^{-\epsilon t} \int_{0}^t \gamma e^{\epsilon s} \rho(\delta, s) \dd s + \gamma \rho(\delta, t)\ .  
\end{equation}
Hence, the boundary condition can be interpreted as a Robin-type boundary conditions with a memory term or a delayed boundary condition
\begin{equation}\label{eq_delayed_bc}
     \pp_x (x (1 - x) \rho) (\delta, t) = - \epsilon e^{-\epsilon t} \int_{0}^t \gamma e^{\epsilon s} \rho(\delta, s) \dd s + \gamma \rho(\delta, t)   
\end{equation}
A similar calculation can be done for $x = 1 - \delta$.

\begin{remark}
If $\epsilon = 0$, the equation becomes
\begin{equation}\label{Eq_Final_LMA_epsilon_0}
  \begin{cases}
    &  \pp_t \rho = \pp_{xx} (x (1 - x) \rho), \quad   x \in (\delta, 1 - \delta) \\
   &  \pp_x (x (1 - x) \rho) |_{x = \delta} = \gamma \rho(\delta, t), \quad    \pp_x (x (1 - x) \rho) |_{x = 1 - \delta}  = - \gamma \rho(1 - \delta, t) \\
  \end{cases}
 \end{equation}
 which can be viewed as a closed equation on $(\delta, 1-\delta)$ with Robin boundary condition. Although the energy-dissipation law (\ref{ED_RK}) with (\ref{Form_G}) is no longer valid, the system can be interpreted as weighted $L^2$-type gradient flow 
\begin{equation}
 \frac{\dd}{\dd t} \left( \int_{\delta}^{1 - \delta} |\pp_x ( x(1 - x) \rho ) |^2\dd x + \gamma \delta(1 - \delta) (|\rho(\delta, t)|^2 + |\rho(1 - \delta, t)|^2)  \right) = - \int_{\delta}^{1-\delta} x(1 - x) |\rho_t|^2 \dd x
\end{equation}
 One can further define
 \begin{equation}
 a'(t) =  \gamma \rho(\delta, t), \quad b'(t) =  \gamma \rho(1 - \delta, t) 
 \end{equation}
 as the fixation dynamics on the boundary.
 \end{remark}

 \begin{remark}
An important feature of the original Kimura equation is the conservation of the first moment, or fixation probability, that is, $ \frac{\dd}{\dd t} \int x \rho \dd x = 0$. For the regularized system (\ref{Eq_Final_LMA}), a direct calculation shows that
\begin{equation}
\begin{aligned}
 &  \frac{\dd}{\dd t} \int_{\delta}^{1 - \delta} x \rho(x, t) \dd x + a \psi(0)  + b \psi(1) \\
 & = \int_{\delta}^{1 - \delta}  - \pp_x (x) \pp_x ( x(1 - x) \rho ) \dd x -  x (\rho u) \Big |_{\delta}^{ 1 - \delta} + a'(t) \psi(0)  + b'(t) \psi(1) \\
 & =  \int_{\delta}^{1 - \delta}  \pp_{xx} (x) ( x(1 - x) \rho ) \dd x - \pp_x (x) (x(1 - x) \rho ) \Big|_{\delta}^{1 - \delta}  + \delta (b'(t) - a'(t))  \\
 & = - \delta (1 - \delta) (\rho(1 - \delta) - \rho(\delta)) + \delta (b'(t) - a'(t)) \\
 & = - \delta (1 - \delta - \gamma) (\rho(1 - \delta) - \rho(\delta))      + \epsilon \delta (a(t) - b(t)) , \\
\end{aligned}
\end{equation}
where $\psi(x) = x$.
Hence, the regularized system (\ref{Eq_Final_LMA}) may not satisfy the conservation of the first moment for given $\delta > 0$ and $\epsilon > 0$. However, the variation will be very small, smaller than the order of $\delta$.
In the case with $\epsilon = 0$, the conservation of the first moment requires $\gamma = 1 - \delta$, which is consistent with the analysis in \cite{chalub2021gradient} for $\delta = 0$.
\end{remark}

\begin{remark}
In the current study, we view $a(t)$ and $b(t)$ as the probability at $x = 0$ and $x = 1$. Alternatively, one can define the probability density $\tilde{\rho}(x, t)$ at $[0, 1]$ from the solution of (\ref{Eq_Final_LMA}) by
\begin{equation}
  \tilde{\rho}(x, t) = 
  \begin{cases}
    & a(t) /  \delta, \quad x \in [0, \delta) \\
    & \rho(x. t), \quad x \in [\delta, 1 - \delta] \\
    & b(t) / \delta, \quad x \in (1 - \delta, 1] \\
    \end{cases}
    \end{equation}
Then the conservation of mass is the same as in (\ref{conservation_of_mass}). The time evolution of the first moment can be computed as follows: 
\begin{equation}
\begin{aligned}
 &  \frac{\dd}{\dd t} \int_{\delta}^{1 - \delta} x \rho(x, t) \dd x + \int_{0}^{\delta} \frac{a}{\delta} x \dd x  +  \int_{1 - \delta}^{1} \frac{b}{\delta} x \dd x \\
 & = \int_{\delta}^{1 - \delta}  - x \pp_x ( x(1 - x) \rho ) \dd x - x (\rho u) \Big |_{\delta}^{ 1 - \delta} + a'(t) \frac{\delta}{2}  + b'(t) (1 - \frac{\delta}{2}) \\
 & =  -  (x(1 - x) \rho ) \Big|_{\delta}^{1 - \delta}  + \frac{\delta}{2} (b'(t) - a'(t))
 = - \delta (1 - \delta - \tfrac{\gamma}{2}) (\rho(1 - \delta) - \rho(\delta))      + \frac{\epsilon \delta}{2} (a(t) - b(t))\ . \\
\end{aligned}
\end{equation}
Again, the regularized system (\ref{Eq_Final_LMA}) may not satisfy the conservation of the first moment for given $\delta > 0$ and $\epsilon > 0$, but the variation in the first moment will be very small. From this point of view, we need to take $\gamma = 2 - 2 \delta$ and $\epsilon = 0$ to guarantee the conservation of the fixation probability.
\end{remark}

Since the goal is to study the regularized model with small $\epsilon$ and $\delta$, in the remainder of this paper, we take $\gamma = 1$ without losing generality.

\subsection{Equilibrium of the regularized model}
 
Next, we perform some formal analysis to the regularized model. At the equilibrium, we have
\begin{equation}\label{ab_inf}
   \rho^{\infty}(\delta) = \epsilon a^{\infty}, \quad \rho^{\infty} (1 - \delta) = \epsilon b^{\infty},
\end{equation}
and $\rho^{\infty}$ satisfies 
\begin{equation}\label{Eq_steady_state}
  \pp_{xx} ( x(1 - x) \rho^{\infty}) = 0, \quad  x \in (\delta, 1 - \delta)\ .
\end{equation}
Also, due to mass conservation, we have
\begin{equation}\label{conservation_mass_eq}
  \int_{\delta}^{1 - \delta} \rho^{\infty} \dd x + a^{\infty} + b^{\infty} = 1 \, 
\end{equation}

It is straightforward to show that for $\delta > 0$, the classical solution to (\ref{Eq_steady_state}) is given by
\begin{equation}
  \rho^{\infty}(x) = \frac{Ax + B}{x (1 - x)}\ , \quad x \in (\delta, 1 - \delta)
\end{equation}
Then according to (\ref{ab_inf}), $A$ and $B$ satisfies
\begin{equation}
  A \delta + B  = \delta (1 - \delta)  \epsilon a^{\infty}, \quad A (1 - \delta) + B = \delta (1 - \delta) \epsilon b^{\infty} \ .
\end{equation}
One can solve $A$ and $B$ in terms of $a^{\infty}$ and $b^{\infty}$, that is,
\begin{equation}\label{value_A_B}
  \begin{aligned}
   & A =  \epsilon \frac{\delta (1 - \delta)}{1 - 2 \delta} (b^{\infty} - a^{\infty})\ , \quad B = \epsilon \frac{\delta (1 - \delta)}{1 - 2 \delta} ( (1 - \delta) a^{\infty} - \delta b^{\infty})\ . \\
  \end{aligned}
\end{equation}
Hence, for fixed $\delta$, by letting $\epsilon \rightarrow 0$, the equilibrium solution $\rho^{\infty}$ goes to $0$.  

By a direct calculation, we have
\begin{equation}
  \begin{aligned}
  \int_{\delta}^{1 - \delta} \rho^{\infty} \dd x  &  = \int_{\delta}^{1 - \delta} \frac{B}{x}  + \frac{A + B}{1 - x}    \dd x  =  B \ln x \Big|_{\delta}^{1 - \delta} - (A + B) \ln (1 - x) |_{\delta}^{1 - \delta } \\
  & 
  = (2 B + A) (\ln (1 - \delta) - \ln \delta) = \epsilon \delta (1 - \delta) (a^{\infty} + b^{\infty}) (\ln (1 - \delta) - \ln \delta)\ .  \\ 
  \end{aligned}
\end{equation}
Using (\ref{conservation_mass_eq}), we have 
\begin{equation}
  a^{\infty} + b^{\infty} = \frac{1}{1 + \epsilon \delta (1 - \delta)  (\ln (1 - \delta) - \ln \delta)  }.
\end{equation}
Therefore, for fixed $\delta$, by taking $\epsilon \rightarrow 0$, we also have $a^{\infty} + b^{\infty}  \rightarrow 1$. Moreover, for fixed $\epsilon$, since $\delta (1 - \delta)  (\ln (1 - \delta) - \ln \delta)  \rightarrow 0$ as $\delta \rightarrow 0$, we also have $a^{\infty} + b^{\infty}  \rightarrow 1$, then the mass conservation indicates that $\rho^{\infty}  \rightarrow 0$.

\section{Existence of Solutions: Operator Approach}
In this section, we prove the existence of global solutions to the regularized Kimura equation (\ref{eq_reg_kimura_intro}) for $\epsilon>0$ and $\delta>0$.
The method we present here is an operator approach introduced by Friedman \cite{friedman1964partial} and has been used for memory-dependent boundary conditions in \cite{hu1996critical,anderson2011global,deng2012blow}. 
Other methods, such as Galerkin methods \cite{knopf2021phase}, can also be used to prove the existence of strong solutions.

The starting point is the regularized Kimura equation, now rewritten into the form
\begin{align}\label{eq: reg Kimura}
    & \partial_t \rho = \partial_x\big( x(1-x)\partial_x \rho\big) + (1-2x)\partial_x \rho -2\rho,\quad \textnormal{in } (\delta,1-\delta)\times(0,T),\\ \label{eq: Kimura ic}
    & \rho(x, 0) =\rho_0(x), \quad \textnormal{on } (\delta,1-\delta),
\end{align}
with the Robin-type boundary conditions with a memory term
\begin{align}\label{eq: Kimura bc}
\begin{split}
     \delta(1-\delta) \partial_n \rho +2\delta \rho &= \epsilon \textnormal{e}^{-\epsilon t}a_0 + \int_0^t  \textnormal{e}^{-\epsilon (t-s)}\rho(\delta,s)\,\textnormal{d}s,\quad \textnormal{on } \delta\times (0,T),\\
      \delta(1-\delta) \partial_n \rho +2\delta \rho &= \epsilon \textnormal{e}^{-\epsilon t}b_0 + \int_0^t  \textnormal{e}^{-\epsilon (t-s)}\rho(1-\delta,s)\,\textnormal{d}s,\quad \textnormal{on } 1-\delta\times (0,T),
\end{split}
\end{align}
where $n$ denotes the outer unit normal, $a_0$ and $b_0$ are initial conditions of $a(t)$ and $b(t)$.

\begin{theorem}
    Assume that the parameters $\delta,\, \epsilon>0$ and let the initial data satisfy $\rho_0\in C(\delta,1-\delta)$ and let $|a_0|,\,|b_0|<\infty$ be bounded.
    Then, for all times $T>0$, there exists a unique classical solution $\rho$ to equations \eqref{eq: reg Kimura}-\eqref{eq: Kimura bc} such that 
    $$\rho\in C^{2,1}((\delta,1-\delta)\times(0,T))\cap C([\delta,1-\delta]\times[0,T]).$$
\end{theorem}

\begin{proof}
The goal of the first step is to show the local existence of classical solutions.
We use a standard approach based on the Green's function formulation. 
For such purposes, we begin by defining $G_R(x, y, t, \tau)$ as the Green’s function for the above equation with the homogeneous Robin boundary condition.
The existence of the Green's function for Robin boundary conditions was shown e.g. in \cite{marras2009continuous,choi2014green}.
We prove the local existence of solutions for the regularized Kimura equation \eqref{eq: reg Kimura}-\eqref{eq: Kimura bc} via a fixed point argument. 
Despite the presence of a memory term and due to the "good" coercive sign of the boundary condition, the steps are nearly identical to well-known results for corresponding localized problems \cite{friedman1964partial}. 
Thus, we do not repeat all the details here. 

Let $\rho\in C([\delta,1-\delta]\times [0,T])$ and define the operator $\mathcal{T}$ via
\begin{align}\label{eq: operator def}
\begin{split}
    \mathcal{T}[\rho](x,t)&:= \int_\delta^{1-\delta} G_R(x,y,t,0) \rho_0(y)\,\textnormal{d}y \\
     &+\int_0^t G_R(x,1-\delta,t,\tau) \bigg(\frac{\epsilon}{\delta(1-\delta)}\int_0^\tau \textnormal{e}^{-\epsilon(\tau-r)}\rho(1-\delta,r)\textnormal{d}r + \frac{\epsilon}{\delta(1-\delta)} \textnormal{e}^{-\epsilon\tau}b_0\bigg)\textnormal{d}\tau\\
     &-\int_0^t G_R(x,\delta,t,\tau) \bigg(\frac{\epsilon}{\delta(1-\delta)}\int_0^\tau \textnormal{e}^{-\epsilon(\tau-r)}\rho(\delta,r)\textnormal{d}r + \frac{\epsilon}{\delta(1-\delta)} \textnormal{e}^{-\epsilon\tau}a_0\bigg)\textnormal{d}\tau.
     \end{split}
\end{align}
Moreover, we define 
\begin{align*}
    M_0&:= \sup_{x\in [\delta,1-\delta]}|\rho_0|,\\
    \nu(t)&:=\sup_{(x,\tau)\in [\delta,1-\delta]\times [0,t]} \int_0^t \int_\delta^{1-\delta} |G_R(x,y,\tau,s)|\textnormal{d}y\textnormal{d}s,\\
    \kappa(t)&:= \sup_{x\in [\delta,1-\delta]\times [0,t]}  \int_0^t |G_R(x,y,\tau,s)| \textnormal{d}s, \quad y\in \{\delta,1-\delta\},
\end{align*}
where the quantities are well-defined by \cite{deng1992influence,hu1996critical} and there exists a constant $C_0>0$ such that $\kappa(t)\leq C_0\sqrt{t}$.
The idea is now to show that the mapping $\mathcal{T}$ is a compact and continuous self-mapping, which then allows for the application of a standard fixed-point argument.
Thus, we estimate
\begin{align*}
   & \sup_{(x,t)\in[\delta,1-\delta]\times [0,T] }|\mathcal{T}[\rho]|\leq \sup_{(x,t)\in[\delta,1-\delta]\times [0,T] }\bigg|  \int_\delta^{1-\delta} G_R(x,y,t,0) \rho_0(y)\,\textnormal{d}y \bigg|\\
    &+ \sup_{(x,t)\in[\delta,1-\delta]\times [0,T] }\bigg|\int_0^t G_R(x,1-\delta,t,\tau) \bigg(\frac{\epsilon}{\delta(1-\delta)}\int_0^\tau \textnormal{e}^{-\epsilon(\tau-r)}\rho(1-\delta,r)\textnormal{d}r + \frac{\epsilon}{\delta(1-\delta)} \textnormal{e}^{-\epsilon\tau}b_0\bigg)\textnormal{d}\tau\bigg|\\
    &+ \sup_{(x,t)\in[\delta,1-\delta]\times [0,T] }\bigg|\int_0^t G_R(x,\delta,t,\tau) \bigg(\frac{\epsilon}{\delta(1-\delta)}\int_0^\tau \textnormal{e}^{-\epsilon(\tau-r)}\rho(\delta,r)\textnormal{d}r + \frac{\epsilon}{\delta(1-\delta)} \textnormal{e}^{-\epsilon\tau}a_0\bigg)\textnormal{d}\tau\bigg|
\end{align*}
Denoting $f_0=\max\{a_0,b_0\}$ we can further estimate that
\begin{align*}
    \sup_{(x,t)\in[\delta,1-\delta]\times [0,T] }|\mathcal{T}[\rho]|\leq M_0 \nu(T)+ C_0\frac{2\epsilon}{\delta(1-\delta)}f_0 \sqrt{T} +M_0C_0\frac{\epsilon}{\delta(1-\delta)}\sqrt{T^3}.
\end{align*}
Now, choosing the time $T^*$ such that $\nu(T^*)+ C_0\frac{2\epsilon}{\delta(1-\delta)}\frac{f_0}{M_0} \sqrt{T^*} +C_0\frac{\epsilon}{\delta(1-\delta)}\sqrt{{T^*}^3}<1$ we have shown that $\mathcal{T}: K\to K$ compactly, where 
\begin{align*}
    K=\{\rho\in C([\delta,1-\delta]\times [0,T]):\|\rho\|_\infty\leq M_0\}.
\end{align*}
To show the continuity of $\mathcal{T}$ we consider two functions $\rho_1,\rho_2\in K$ with initial data $\rho_{0,1},\rho_{0,2}$ and boundary data $a_{0,1},a_{0,2},b_{0,1},b_{0,2}$.
Then, it follows from the linearity of \eqref{eq: operator def} that $|\mathcal{T}[\rho_1-\rho_2]|\to 0$ when $|\rho_{0,1}-\rho_{0,2}|\to 0$ and $|a_{0,1}-a_{0,2}|,|b_{0,1}-b_{0,2}|\to 0$. 
Therefore the operator $\mathcal{T}$ has a fixed point and this guarantees the existence of local-in-time classical solutions.

Next, to show the positivity of solutions we apply a comparison principle.
Consider $w(x,t)\in C^{2,1}((\delta,1-\delta)\times(0,T))\cap C([\delta,1-\delta]\times[0,T]) $ such that
\begin{align*}
    \partial_t w&\geq \partial_{x,x}(x(1-x)w),\quad \textnormal{in } (\delta,1-\delta)\times (0,T),\\
    \partial_n w&\geq -2\delta w +\epsilon\textnormal{e}^{-\epsilon t}a_0 +\int_0^t \textnormal{e}^{-\epsilon(t-s)}w(\delta,s)\textnormal{d}s,\quad \textnormal{on } \{\delta\}\times(0,T),\\
     \partial_n w&\geq -2\delta w +\epsilon\textnormal{e}^{-\epsilon t}b_0 +\int_0^t \textnormal{e}^{-\epsilon(t-s)}w(1-\delta,s)\textnormal{d}s,\quad \textnormal{on } \{1-\delta\}\times(0,T),\\
     w(0)&=w_0\geq 0,\quad \textnormal{on } [\delta,1-\delta].
\end{align*}
For a positive smooth function $\xi$, satisfying $\delta(1-\delta)\partial_n \xi \geq \alpha \xi $ for $x=\delta,1-\delta$ we define a new function $W$ by
\begin{align*}
    w(x,t)=\textnormal{e}^{\lambda t}\xi(x)W(x,t),
\end{align*}
where the constants $\alpha,\lambda$ are chosen such that
\begin{align*}
    \alpha&> -2+\epsilon T, \quad  \lambda>\frac{|\Delta \xi|}{\xi}+2\frac{|\partial_x \xi|}{\xi} -2.
\end{align*}
Then, we obtain
\begin{align*}
    \partial_t W&\geq  x(1-x)\Delta W+\big(x(1-x)\frac{2\partial_x \xi }{\xi}+2(1-2x)\big)\partial_x W\\
    &\quad + \big( x(1-x)\frac{\Delta\xi}{\xi}+ 2(1-2x)\frac{\partial_x \xi}{\xi}-\lambda-2\big)W,
\end{align*}
which implies, using the classical maximum principle, that there is no negative minimum of W in $(\delta,1-\delta)\times [0,T]$.
Now, we assume that there exists a $W(x_0,t_0)=\min W<0$ for some $x_0\in \{\delta,1-\delta\},t_0\in [0,T]$.
It follows that
\begin{align*}
    \delta(1-\delta) \partial_n W\geq - W(x_0,t_0) \bigg[ \delta(1-\delta)\partial_n \xi +2 \xi -\epsilon\xi \textnormal{e}^{-\lambda t}\int_0^t \textnormal{e}^{-\epsilon(t-s)+\lambda s}\textnormal{d}s\bigg]+ \epsilon\textnormal{e}^{-\epsilon t}f_0>0
\end{align*}
which is a contradiction. Hence, $W\geq 0$ in $[\delta,1-\delta]\times [0,T]$, which implies that $w\geq 0$.

By the same reasoning we obtain a comparison principle for the equation.
Let $\overline \rho\geq c>0$ be a supersolution  and $\underline \rho\geq 0$ be a subsolution of the regularized Kimura equation.
Then setting $w= \overline \rho-\underline\rho $ we can apply the above result and obtain that $\overline \rho\geq \underline\rho$ provided $\overline \rho_0\geq \underline\rho_0$.

For the  existence of global solutions we have to show that we can extend the local solution for all times $T>0$, i.e. we have to show that $\rho(x,t)\leq C(T)<\infty$ for $(x,t)\in [\delta,1-\delta]\times [0,T]$. 
To do this we seek a global supersolution of the regularized Kimura equation.
First, we know that there exists a function $\phi(x)\in C^2(\delta,1-\delta)$ such that $0< \phi(x)\leq 1$ in $(\delta,1-\delta)$, $\partial_n \phi\geq 1$ for $x\in \{\delta,1-\delta\}$ \cite{deng2000solutions}.
Then, we define $\overline{\rho}$ as 
\begin{align*}
    \overline{\rho}= M\textnormal{e}^{\lambda t+\phi},
\end{align*}
where 
$M=\max\big\{\frac{\epsilon f_0}{2\delta-\epsilon},\|\rho_0\|_\infty\big\}$, $\lambda=\max\{ m_1^2+2m_1+m_2-2,1\}$ and $m_1= \sup |\partial_x \phi|$, $m_2=\sup|\Delta\phi|$.
It follows that $\overline{\rho}$ satisfies 
\begin{align*}
 \partial_t \overline{\rho}&\geq \partial_{x,x} \big(x(1-x)\overline{\rho}\big),\quad \textnormal{in } (\delta,1-\delta)\times(0,T),\\
 \partial_n \overline{\rho}&\geq -2\delta \overline{\rho}+ \epsilon\textnormal{e}^{-\epsilon t}f_0+\int_0^t \textnormal{e}^{-\epsilon(t-s) \overline{\rho}}\textnormal{d}s,\quad \textnormal{on } \{\delta,1-\delta\}\times (0,T),\\ 
 \overline{\rho}(x,0)&\geq \rho_0, \quad \textnormal{on } [\delta,1-\delta]\times \{0\}.
\end{align*}
Hence, $\overline{\rho}$ is the desired bounded supersolution, which implies the existence of a global solution.
\end{proof}

\section{Numerics}
In this section, we perform a numerical study to the regularized Kimura equation, considering various initial conditions as well as different values of $\delta$ and $\epsilon$.  
These numerical results demonstrate that the regularized equation can capture the key properties of the original Kimura equation.


\subsection{Numerical schemes}
We first propose a numerical scheme for the regularized Kimura equation \eqref{Eq_Final_LMA} based on a finite volume approach. For a given $\delta$, we can divide $(\delta, 1 - \delta)$ into $N$ subintervals of size $h = (1 - 2 \delta) / N$. Let $x_i = \delta + ih$ ($i = 0, 1, 2 \ldots, N$), and define
\begin{equation}
   \rho_0 = \frac{2}{h} \int_{\delta}^{x_{1} + h/2} \rho \dd x,  \quad  \rho_i = \frac{1}{h} \int_{x_i - h/2}^{x_i + h/2} \rho \dd x, \quad i = 1, \ldots N - 1, \quad \rho_N = \frac{2}{h} \int_{x_N - h/2}^{x_N} \rho \dd \x.
\end{equation}
Due to the mass conservation $ \frac{\dd}{\dd t} \left(\int_{\delta}^{1 -\delta} \rho \dd x + a(t) + b(t) \right) = 0$, we have 
\begin{equation}\label{Dis_MC}
  a(t) + b(t) +  h \left( \frac{1}{2} \rho_0(t)  + \sum_{i = 1}^{N-1} \rho_i(t) + \frac{1}{2} \rho_N(t) \right) = 1, \quad \forall t
\end{equation}
in the semi-discrete sense.  Let ${\bm p}(t) = (a, \rho_0, \ldots, \rho_N, b)^{\rm T}$, and ${\bm w} = (1, h/2, h, \ldots, h, h/2, 1) \in \mathbb{R}^{N+3}$. 
It is convenient to define the inner product in the discrete space
\begin{equation}
 \langle f, g \rangle_{w} = \sum_{i=1}^{N+3} w_i f_i g_i.
\end{equation}
Then the mass conservation (\ref{Dis_MC}) can be written as $\langle {\bm p}, {\bf 1} \rangle = 1$. 
We also define the discrete free energy as
\begin{equation}
  \begin{aligned}
  \mathcal{F}_h({\bm p})  & =  h  \left( \frac{1}{2}\rho_0 \ln (\delta (1 - \delta) \rho_0) +  \sum_{i=1}^{N+1} \rho_i \ln ( x_i (1 - x_i) \rho_i) +  \frac{1}{2} \rho_N \ln (\delta (1 - \delta) \rho_N) \right) \\ 
  & + a \ln (\epsilon \delta (1 - \delta) a) + b \ln (\epsilon \delta (1 - \delta) b) = \langle {\bm \mu}, {\bm p} \rangle_{w}, \\
  \end{aligned}
\end{equation}
where ${\bm \mu} = (\ln (\epsilon \delta (1 - \delta) a), \ln (x_0 (1 - x_0) \rho)_0), \ldots \ln (x_N (1 - x_N) \rho)_N)  ,\ln (\epsilon \delta (1 - \delta) b))^{\rm T}$ is the discrete chemical potential.

The semi-discrete system associated with (\ref{Eq_Final_LMA}) can be written as 
\begin{equation}\label{Eq_Final_semi}
  \begin{cases} 
  & \frac{\dd a}{\dd t}  =   (\rho_0 - \epsilon a) \\
  & \frac{\dd \rho_0}{\dd t} = \frac{2}{h}  ( - (\rho u)_{1/2}  + (\epsilon a - \rho_0)  \\ 
  & \frac{\dd \rho_i}{\dd t}  = \frac{1}{h}  (  - (\rho u)_{i + 1/2}  + (\rho u)_{i-1/2}), \quad i = 1, \ldots N -1  \\ 
  & \frac{\dd \rho_N}{\dd t} = \frac{2}{h} (  (\epsilon b - \rho_N)  + (\rho u)_{N - 1/2}) \\
  & \frac{\dd b}{\dd t}  = (\rho_N  - \epsilon b) \ , \\
  \end{cases}
\end{equation}
where
\begin{equation}
  (\rho u)_{i + 1/2}  = - \frac{1}{h} (x_{i+1} (1 - x_{i+1}) \rho_{i+1} - x_{i} (1 - x_{i}) \rho_{i} ) \ . 
\end{equation}
it is clear that (\ref{Eq_Final_semi}) is a linear system of ${\bm p}$, which can be denoted by
\begin{equation}
  \frac{\dd {\bm p}}{\dd t} = L_h {\bm p},  \quad L_h \in \mathbb{R}^{(N+3) \times (N+3)}.
\end{equation}
Moreover, since $\langle \frac{\dd p}{\dd t}, 1 \rangle_{w} = \langle L_h {\bm p}, 1 \rangle_w$,  the semi-discrete equation (\ref{Eq_Final_semi}) preserves the discrete mass conservation (\ref{Dis_MC}).

Next, we introduce a temporal discretization to (\ref{Eq_Final_semi}). Since we are interested in investigating the behavior of $\epsilon$ and $\delta$, a high-order temporal discretization is required. 
In the current study, we adopt a second-order Crank--Nicolson scheme, given by
\begin{equation}\label{CN}
  \frac{{\bm p}^{n+1} - {\bm p}^n}{\tau} = \frac{1}{2} (L_h {\bm p}^{n} + L_h {\bm p}^{n+1}),
\end{equation}
for the temporal discretization. Since the scheme is linear, we can solve $p^{n+1}$ directly, given by 
\begin{equation}
  {\bm p}^{n+1} =  \left( \frac{1}{\tau} {\sf I} - \frac{1}{2} L_h \right)^{-1}  \left( \frac{1}{\tau} {\sf} + \frac{1}{2}L_h \right) {\bm p}^n\ ,
\end{equation}
provided the matrix $\frac{1}{\tau} {\sf I} - \frac{1}{2} L_h$ is  invertible.
Although it is not straightforward to prove, numerical simulations show that the numerical scheme (\ref{CN}) is positive-preserving and energy stable. We'll address this in the future work.

\subsection{Numerical results}
In this subsection, we present some numerical results for the regularized Kimura equation with different initial conditions and and various values of $\epsilon$ and $\delta$.

\begin{figure}[!h]
\begin{overpic}[width = 0.48 \linewidth]{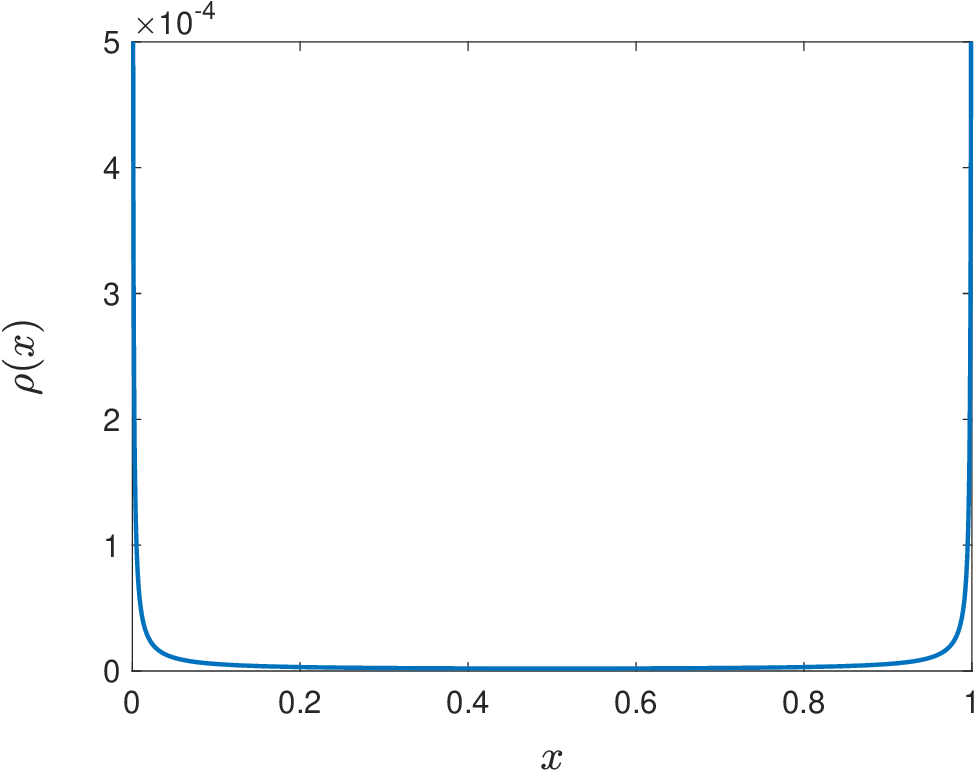}
\put(0, 70){(a)}
\end{overpic}
\hfill
\begin{overpic}[width =  0.44 \linewidth]{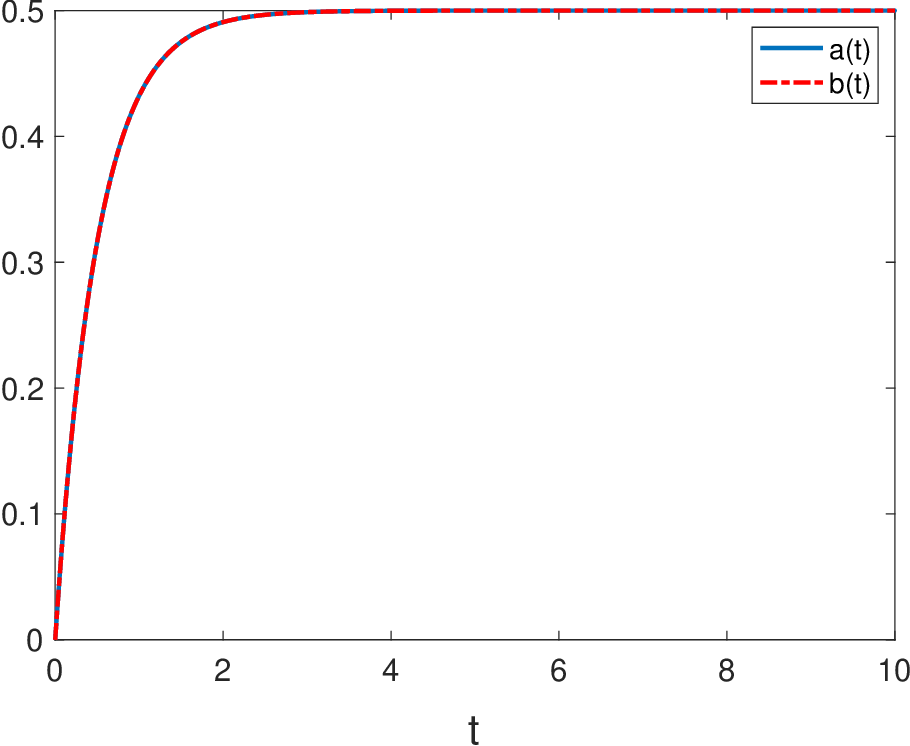}
\put(-9, 78){(b)}
\end{overpic}

\vspace{1em}
\begin{overpic}[width = 0.48 \linewidth]{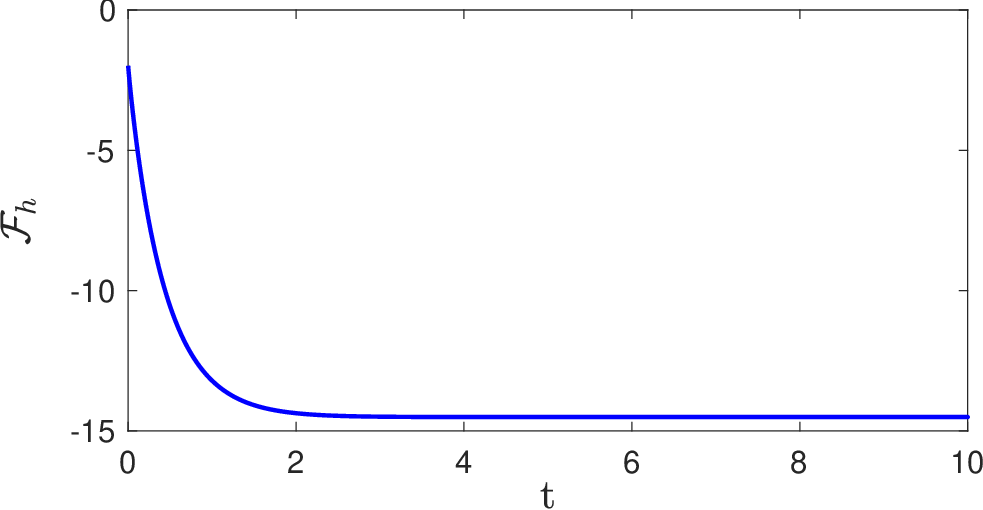}
\put(0, 48){(c)}
\end{overpic}
\hfill
\begin{overpic}[width = 0.48 \linewidth]{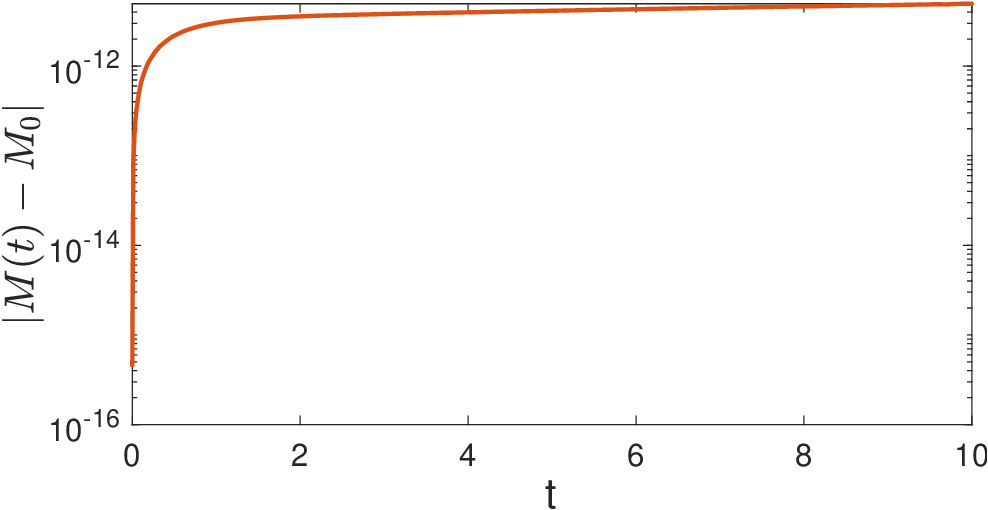}
\put(-2, 48){(d)}
\end{overpic}
  \caption{(a) Numerical solutions at $t = 10$ for $\delta = 0.001$ with $\epsilon = 0.001$ with the uniform initial condition (\ref{IC1}). (b) Numerical solution for $a(t)$ and $b(t)$.  (c) Evolution of the discrete energy. (d) Time evolution of the absolute difference between the first moment $M(t)$ and its initial value $M_0$. }\label{Fig1}
  \end{figure}

\noindent {\bf Uniform boundary condition}: We first consider a uniform initial condition, given by
\begin{equation}\label{IC1}
  \rho(x, 0) = 1 / (1 - 2 \delta), \quad x \in (\delta, 1 - \delta), \quad a(0) = b(0) = 0. 
 \end{equation}
Fig. \ref{Fig1}(a) shows the numerical solution $\rho(x, t)$ for $\delta = 10^{-3}$ and $\epsilon = 10^{-3}$ at $t = 10$, with $h = 10^{-4}$ and $\tau = 10^{-4}$. The time evolution of $a(t)$ and $b(t)$ is shown in Fig. \ref{Fig1} (b). Due to the symmetry of the initial condition, the dynamics of $a(t)$ and $b(t)$ is exactly the same.
It is clear that nearly all of the mass is moving to the boundary. Figures \ref{Fig1}(c) and \ref{Fig1}(d) show the evolution of the discrete energy and the absolute difference of the first moment $M(t)$ from its initial value $M_0$, respectively.  It can be noticed that the discrete energy decreases with respect to time, while the first moment is conserved numerically.
We need to mention that, due to the symmetry of the initial condition, we have $a'(t) = b'(t)$ and $\rho(1 - \delta, t) = \rho(\delta, t)$ for all $t$, so the first moment is also conserved in theory.

\begin{figure}[!h]
  \includegraphics[width = 0.48 \linewidth]{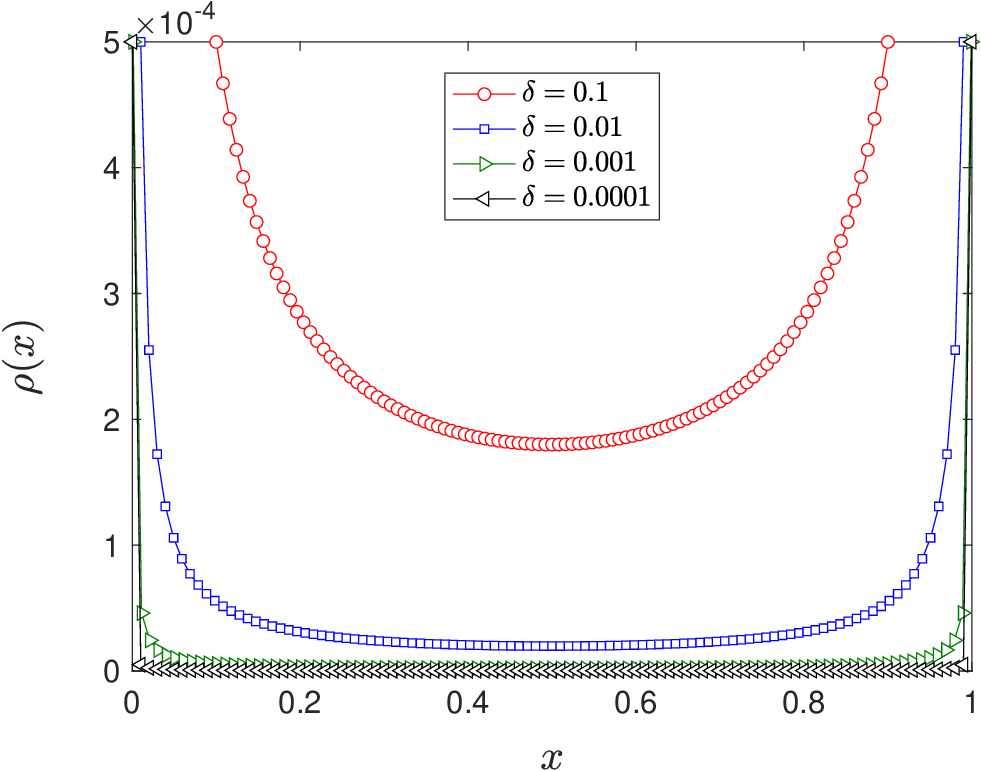}
  \hfill
  \includegraphics[width = 0.48 \linewidth]{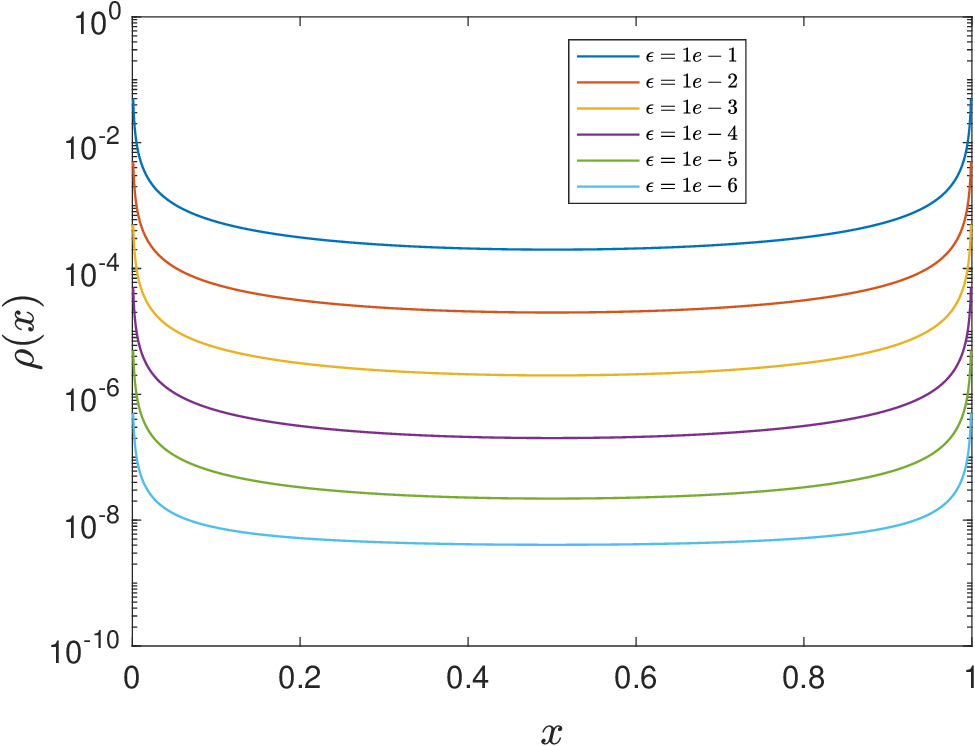}
  
  \caption{(a) Numerical solution in the bulk ($\rho(x, t)$) at $t = 10$ for fixed $\epsilon = 10^{-3}$ and varying $\delta$. (b) Numerical solution in the bulk ($\rho(x, t)$) at $t = 10$ for fixed $\delta = 10^{-3}$ and different $\epsilon$. $\Delta t = 10^{-4}$ and $h = (1 - 2\delta) / 10^4$.}\label{Fig_delta_epsilon}
  \end{figure}

Next, we conducted numerical experiments to investigate the influence of $\epsilon$ and $\delta$ on the numerical solution. We take $h = 10^{-4}$ and $\tau = 10^{-4}$. Figure \ref{Fig_delta_epsilon}(a) shows the numerical solutions $\rho(x, t)$ at $t = 10$ for various values of $\delta$, while keeping $\epsilon$ fixed at $10^{-3}$. Figure \ref{Fig_delta_epsilon}(b) shows the numerical solutions at $t = 10$ for different values of $\epsilon$ with $\delta$ fixed at $10^{-3}$. Since the scheme is second-order accurate in both time and space, the numerical error is expected to be smaller than $O(10^{-8})$, which has a negligible effect on the plots of $\rho(x, t)$ in Fig.~\ref{Fig_delta_epsilon}. The obtained results clearly demonstrate that the value at $x = \delta$ and $1-\delta$ is determined by $\epsilon$. As $\epsilon$ approaches zero, the bulk solution $\rho(x, t)$ at $t = 10$ tends to zero. The simulation results suggest that if $\epsilon \rightarrow 0$, $\rho(x, t) \rightarrow 0$ when $t \rightarrow 0$ for $x \in (\delta, 1 - \delta)$.

\noindent {\bf Non-uniform initial condition}:
Next, we consider a non-uniform initial condition, given by
\begin{equation}\label{IC2}
\rho(x) = 
\begin{cases}
  &  0.5 / ( 1 - 2 \delta), \quad \delta < x < 0.5 \\
  & 1.5 /  ( 1 - 2 \delta), \quad 0.5 <= x < 1 - \delta \ , \\
\end{cases} \quad a(0) = b(0) = 0.
\end{equation}

\begin{figure}[!h]
\begin{overpic}[width = 0.48 \linewidth]{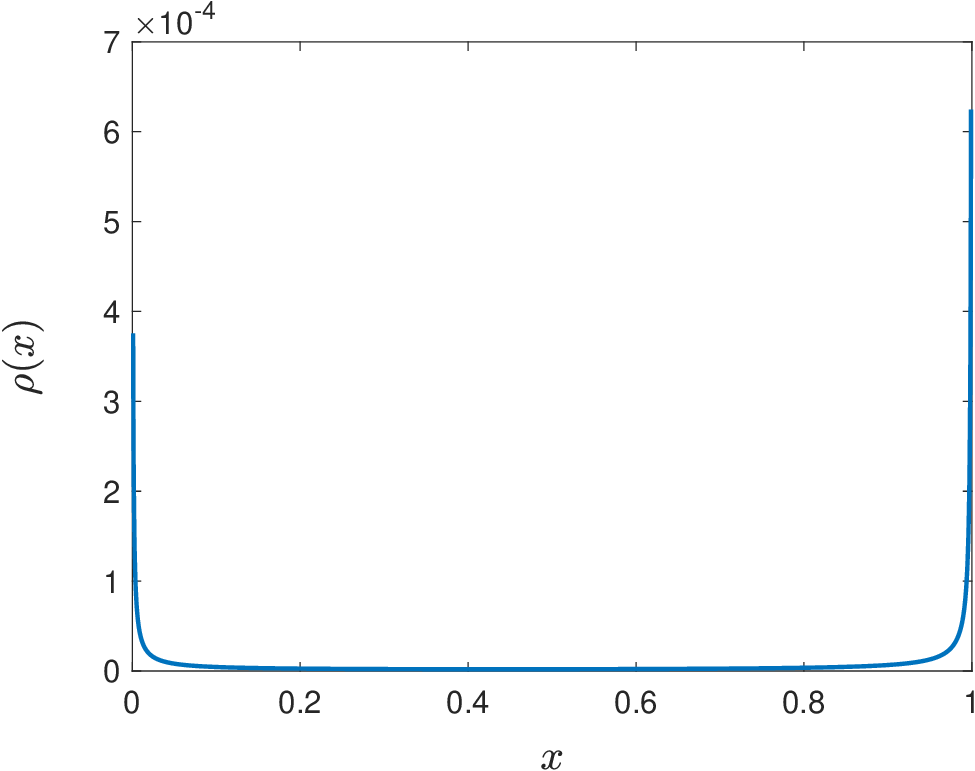}
\put(0, 70){(a)}
\end{overpic}
\hfill
\begin{overpic}[width = 0.44 \linewidth]{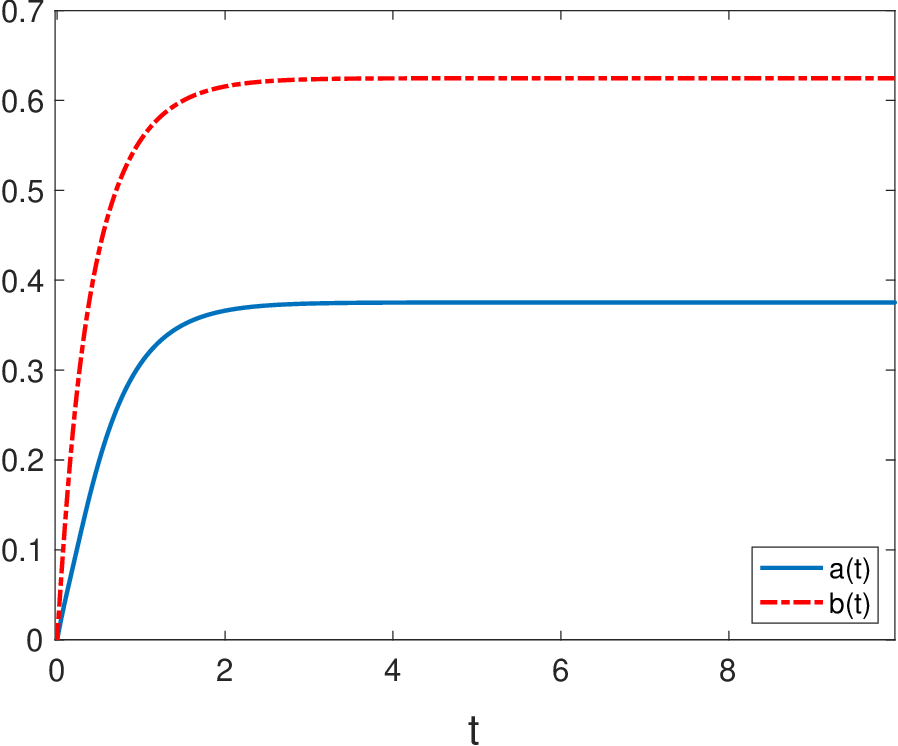}
\put(-9, 78){(b)}
\end{overpic}

\vspace{1em}
\begin{overpic}[width = 0.48 \linewidth]{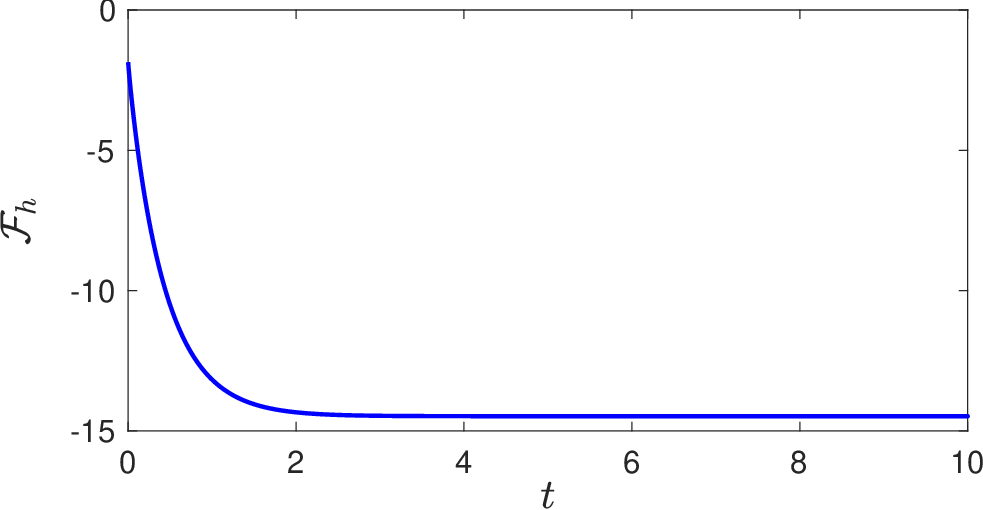}
\put(0, 48){(c)}
\end{overpic}
\hfill
\begin{overpic}[width = 0.48 \linewidth]{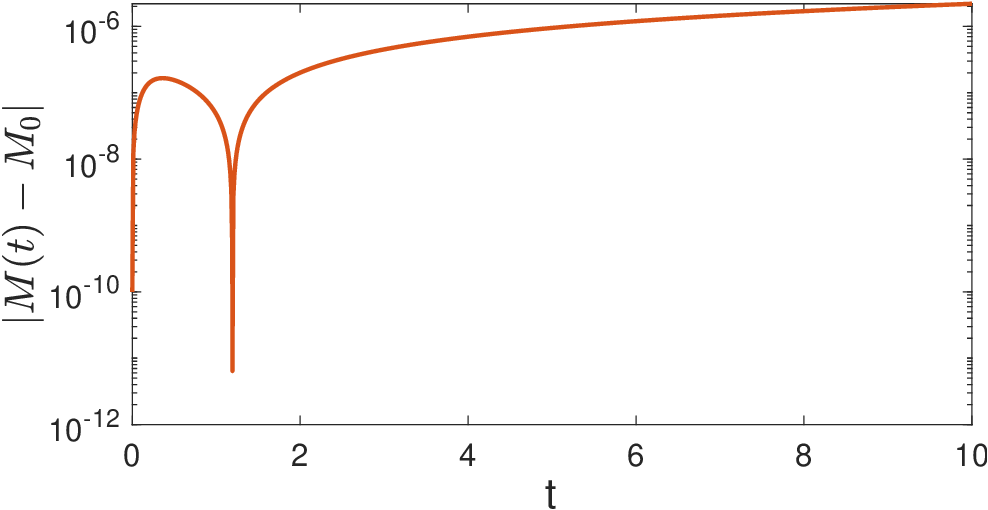}
\put(-2, 48){(d)}
\end{overpic}

\caption{ (a) Numerical solution $\rho(x, t)$ at $t = 10$ for $\delta = 10^{-3}$ and $\epsilon = 10^{-3}$, using the initial condition (\ref{IC2}).  (b) Numerical solution for $a(t)$ and $b(t)$.  (c) Evolution of the discrete energy. (d) Time evolution of the absolute difference between the first moment $M(t)$ and its initial value $M_0$.}\label{Fig_nonuniform}
\end{figure}

Fig. \ref{Fig_nonuniform}(a) shows the numerical solution $\rho(x, t)$ for $\delta = 0.001$ and $\epsilon = 0.001$ at $t = 10$, with $h = 10^{-4}$ and $\tau = 10^{-4}$. The time evolution of $a(t)$ and $b(t)$ is shown in Fig. \ref{Fig_nonuniform} (b). The evolution of discrete energy and the first moment are shown in Fig. \ref{Fig_nonuniform} (c) and (d). In this case, we can also observe the energy stability, and the first moment is almost a constant.

\begin{figure}[!h]
\begin{overpic}[width = 0.48 \linewidth]{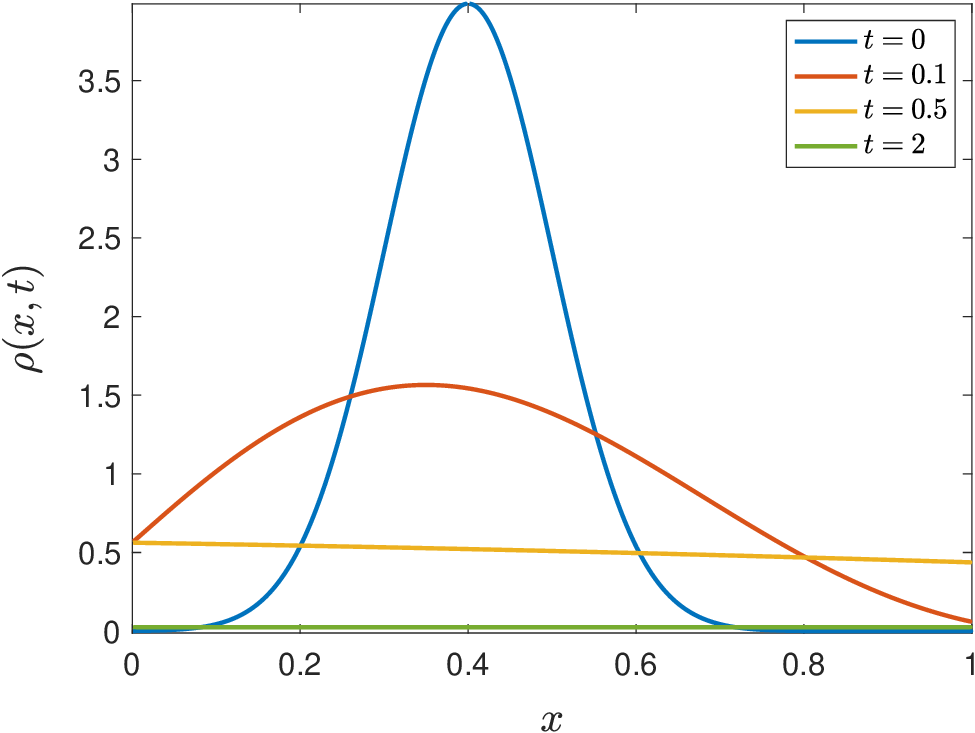}
\put(0, 70){(a)}
\end{overpic}
\hfill
\begin{overpic}[width = 0.44 \linewidth]{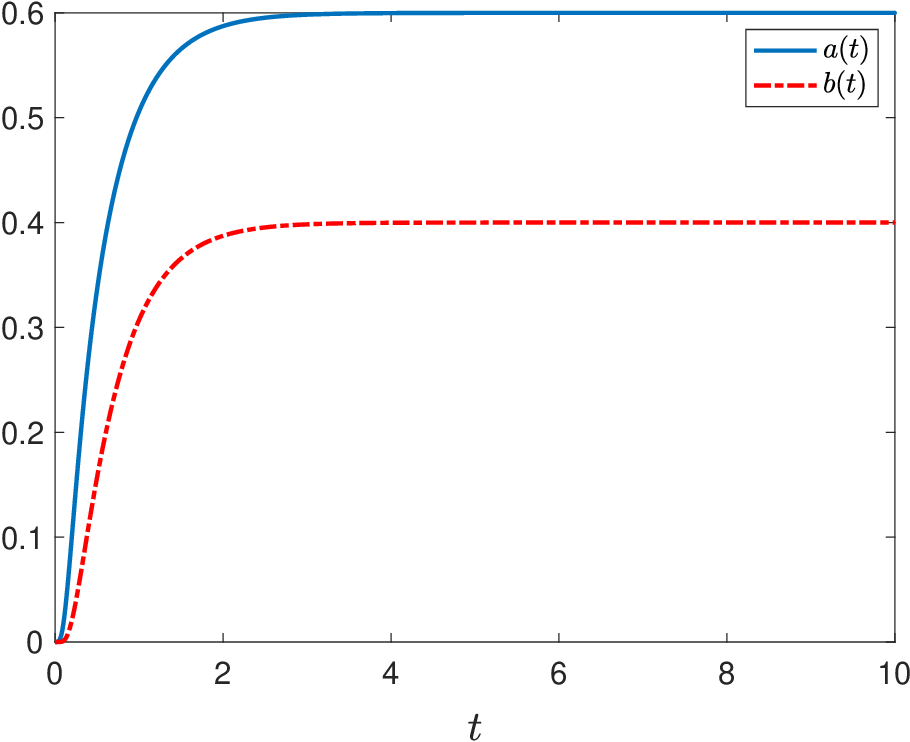}
\put(-9, 78){(b)}
\end{overpic}

\vspace{1em}
\begin{overpic}[width = 0.48 \linewidth]{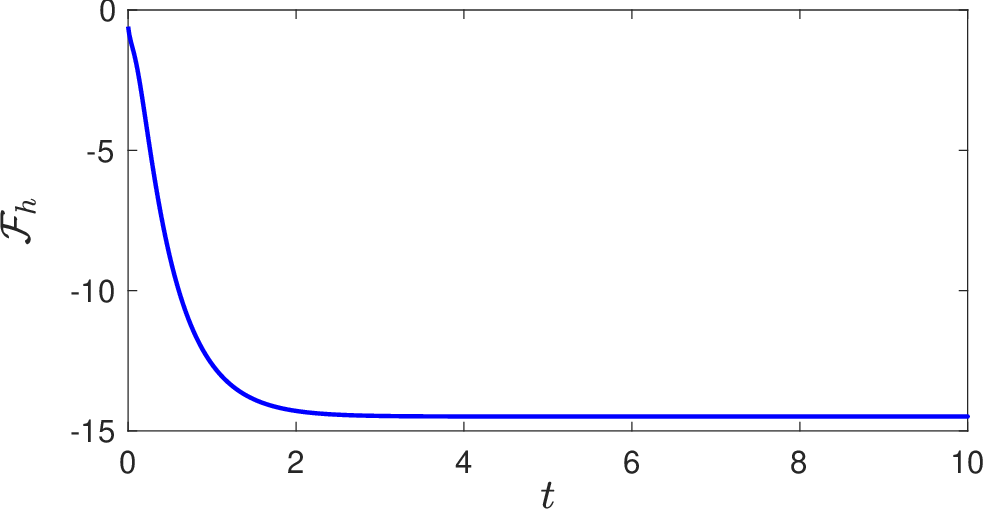}
\put(0, 48){(c)}
\end{overpic}
\hfill
\begin{overpic}[width = 0.48 \linewidth]{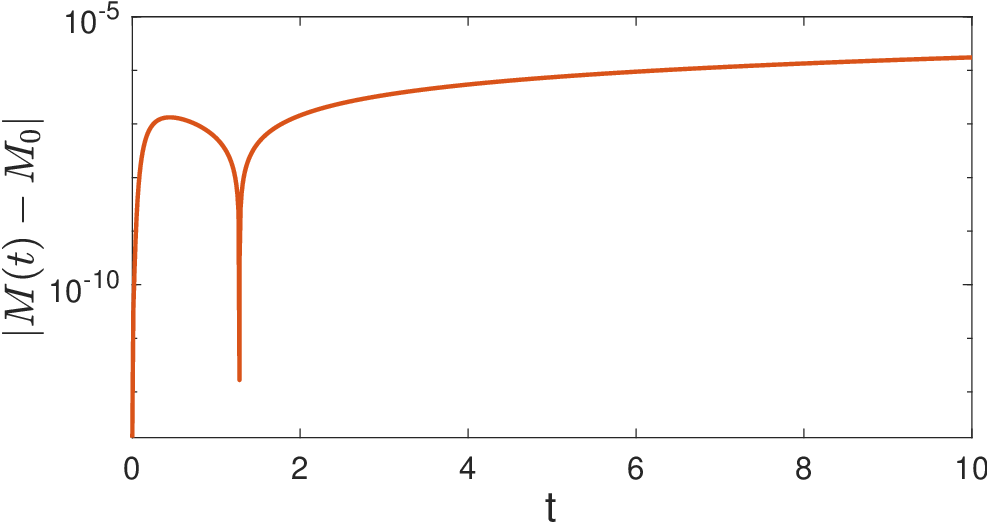}
\put(-3, 48){(d)}
\end{overpic}
\caption{ (a) Numerical solution $\rho(x, t)$ at different $t$ for $\epsilon  = 10^{-3}$ and $\delta = 10^{-3}$, using a Gaussian distribution (\ref{Gaussian_ini}) as the initial condition.(b) Numerical solution for $a(t)$ and $b(t)$. (c) Evolution of the discrete energy. (d) Time evolution of the absolute difference between the first moment $M(t)$ and its initial value $M_0$.} \label{Fig_Gaussian}
\end{figure}

\noindent {\bf Gaussian initial distribution}:
Last, we take the initial condition as a Gaussian distribution 
\begin{equation}\label{Gaussian_ini}
\rho_0(x) = \frac{1}{\sqrt{2 \pi} \sigma} \exp( - \frac{(x - x_0)^2}{2 \sigma^2}), \quad a(0) = b(0) = 0,
\end{equation}
with $\sigma = 0.1$ and $x_0 = 0.4$. Fig. \ref{Fig_Gaussian} shows the numerical results for $\delta = 10^{-3}$ and $\epsilon = 10^{-3}$. Fig. \ref{Fig_Gaussian}(a) shows $\rho(x, t)$ for $t = 0$, $0.1$, $0.5$ and $2$, respectively, while Fig. \ref{Fig_Gaussian}(b) shows the evolution of $a(t)$ and $b(t)$ with respect to $t$.
The evolution of the discrete energy and the first moment are shown in Fig. \ref{Fig_nonuniform} (c) and (d). Again, the numerical scheme is energy stable and first moment is almost conserved.

\section{Conclusion remark}
We proposed a new continuum model for a random genetic drift problem by incorporating a dynamic boundary condition approach. The dynamic boundary condition compensates for singularities on the boundary in the original Kimura equation. We have demonstrated the existence and uniqueness of a strong solution for the regularized system. Finally, we present some numerical results for the regularized model, which indicate that the model can capture the main features of the original model. As future work, we will further study the long-term behavior of the new model.  Additionally, we plan to extend current approach to multi-alleles genetic drift problem.

\section*{Acknowledgements}
C. L. is partially supported by NSF grants DMS-1950868, DMS-2118181 and DMS-2410742. J.-E. S. would like to thank C.L and the Illinois Institute of Technology for the invitation to a research visit. Moreover, J.-E.S. acknowledges the support of the DFG under grant No. 456754695. Y. W is partially supported by NSF DMS-2153029 and DMS-2410740.

\bibliographystyle{siam}
\bibliography{Ref}

\end{document}